\documentclass[fleqn]{mat01}
\usepackage{times,mathtimy,amssymb,latexsym}
\begin{document}

\setcounter{page}{61} \firstpage{61}

\newtheorem{theore}{Theorem}
\renewcommand\thetheore{\arabic{section}.\arabic{theore}}
\newtheorem{theor}[theore]{\bf Theorem}
\newtheorem{lem}[theore]{\it Lemma}
\newtheorem{propo}[theore]{\rm PROPOSITION}
\newtheorem{coro}[theore]{\rm COROLLARY}
\newtheorem{definit}[theore]{\rm DEFINITION}
\newtheorem{probl}[theore]{\it Problem}
\newtheorem{exampl}[theore]{\it Example}
\newtheorem{remar}[theore]{\it Remark}

\newcommand{\D}{\mathcal{D}}
\newcommand{\N}{\mathcal{N}}
\newcommand{\Na}{\mathbb{N}}
\newcommand{\R}{{\mathbb R}}
\newcommand{\Z}{{\mathbb Z}}
\newcommand{\V}{{\mathbb V}}
\newcommand{\U}{{\mathcal U}}
\newcommand{\B}{{\mathcal B}}
\newcommand{\Law}{\mathcal{L}}
\newcommand{\M}{{\widetilde M}}
\newcommand{\VPN}{(V,\nu, \tau,\tau^*)}
\newcommand{\sgn}{\mathop{\rm sgn}\nolimits}
\newcommand{\Int}{mathop{\rm Int}\nolimts}
\newcommand{\Ran}{{\rm Ran}\,}

\newcommand{\norm}[1]{\mid\hspace{-.06cm}\mid \hspace{-.1cm} #1
\hspace{-.1cm}\mid\hspace{-.06cm}\mid}

\newcommand{\h}{\hat \;}
\newcommand{\ie}{i.e.}

\newcommand{\dfe}{d.f.}

\newcommand{\ip}[2]{(#1, #2)}


\title{Quotient probabilistic normed spaces and completeness
results}

\markboth{Bernardo Lafuerza-Guill\'en et al}{Quotient
probabilistic normed spaces}

\author{BERNARDO LAFUERZA-GUILL\'EN$^1$, DONAL O'REGAN$^2$ and
REZA SAADATI$^{3}$}

\address{$^{1}$Departamento de Estad\'\i stica y Matem\'atica
Aplicada, Universidad de Almer\'\i a, 04120~Almer\'\i a, Spain\\
\noindent $^{2}$Department of Mathematics,
National University of Ireland, Galway, Ireland\\
\noindent $^{3}$Faculty of Sciences, University of Shomal, Amol,
Iran\\
\noindent E-mail: blafuerz@ual.es; donal.oregan@nuigalway.ie;
rsaadati@eml.cc}

\volume{117}

\mon{February}

\parts{1}

\pubyear{2007}

\Date{MS received 8 December 2005; revised 30 March 2006}

\begin{abstract}
We introduce the concept of quotient in PN spaces and give some
examples. We prove some theorems with regard to the completeness
of a quotient.
\end{abstract}


\maketitle

\section{Introduction}

In the literature devoted to the theory of probabilistic normed
spaces (PN spaces, briefly), topological and completeness
questions, boundedness and compactness concepts
\cite{LRS1,LRS3,sa}, linear operators, probabilistic norms for
linear operators \cite{LRS4}, product spaces \cite{L} and fixed
point theorems have been studied by various authors. However
quotient spaces of PN spaces have never been considered. This note
is a first attempt to fill this gap.

The present paper is organized as follows. In \S\ref{S:def} all
necessary preliminaries are recalled and notation is established.
In \S\ref{S:quot}, the quotient space of a PN space with respect
to one of its subspaces is introduced and its properties are
studied. Finally, in \S\ref{S:compl}, we investigate the
completeness relationship among the PN spaces considered.

\section{Definitions and preliminaries}\label{S:def}

\setcounter{theore}{0}

In the sequel, the space of all probability distribution functions
(briefly, d.f.'s) is $ \Delta^{+} =\{ F\hbox{:}\ {\bf
R}\cup\{-\infty,+\infty\}\longrightarrow [0,1]\hbox{\rm :}\ F$ is
left-continuous and non-decreasing on $\bf{R}$,  $F(0)=0$ and
$F(+\infty)=1\}$ and the subset $D^{+} \subseteq \Delta^{+}$ is
the set $ D^{+}=\{F\in \Delta^{+}\hbox{:}\ l^{-}F(+\infty)=1\}$.
Here $l^-f(x)$ denotes the left limit of the function $f$ at the
point $x$, $l^- f(x)=\lim_{t\to x^-}f(t)$. The space $\Delta^{+}$
is partially ordered by the usual point-wise ordering of
functions, i.e., $F\leq G$ if and only if $F(x)\leq G(x)$ for all
$x$ in $\bf{R}$. The maximal element for $\Delta^{+}$ in this
order is the d.f. given by
\begin{equation*}
\varepsilon_{0} = \begin{cases}
0, &\text{if} \; x\leq 0,\\[.2pc]
1, &\text{if} \; x>0.
\end{cases}
\end{equation*}\newpage

\noindent Also the minimal element for $\Delta^{+}$ in this order
is the d.f. given by
\begin{equation*}
\varepsilon_{\infty} = \begin{cases}
0, &\text{if} \; x\leq \infty,\\[.2pc]
1, &\text{if} \; x=\infty.
\end{cases}
\end{equation*}
We assume that $\Delta$ is metrized by the Sibley metric $d_S$,
which is the modified L\'{e}vy metric \cite{schw-sk1,Sibley}. If
$F$ and $G$ are \dfe 's  and $h$ is in $(0,1]$, let $(F,G;h)$
denote the condition
\begin{equation*}
F(x-h)-h\leq G(x)\leq F(x+h)+h
\end{equation*}
for all $x$ in $(-1/h,1/h)$. Then the modified L\'{e}vy metric
(Sibley metric) is defined by
\begin{equation*}
d_{S}(F,G):=\inf\{h>0\hbox{\rm :}\  \text{both}\ (F,G;h) \
\text{and} \ (G,F;h) \ \text{hold}\}.
\end{equation*}
For any $F$ in $\Delta^{+}$,
\begin{align*}
d_{S}(F,\varepsilon_{0}) &= \inf\{h>0\hbox{\rm :}\
(F,\varepsilon_{0};h)\ \text{holds}\}\\[.5pc]
&= \inf\{h>0\hbox{\rm :}\ F(h^+)>1-h\},
\end{align*}
and for any $t>0$,
\begin{equation*}
F(t)>1-t\Longleftrightarrow d_{S}(F,\varepsilon_{0})<t.
\end{equation*}
It follows that, for every $F,G$ in $\Delta^{+}$,
\begin{equation*}
F\leq G \Longrightarrow d_{S}(G,\varepsilon_{0})\leq
d_{S}(F,\varepsilon_{0}).
\end{equation*}

A sequence $(F_n)$ of \dfe 's  converges weakly to a d.f. $F$ if
and only if the sequence $(F_n(x))$ converges to $F(x)$ at each
continuity point $x$ of $F$. For the proof of the next theorem see
Theorem~4.2.5 of \cite{schw-sk1}.

\begin{theor}[\!] 
Let $(F_n)$ be a sequence of functions in $\Delta,$ and let $F$ be
in $\Delta$. Then $F_n\rightarrow F$ weakly if and only if
$d_S(F_n,F)\rightarrow 0$.
\end{theor}

\begin{definit}$\left.\right.$\vspace{.5pc}  

\noindent {\rm A \emph{triangular norm} $T$ (briefly, a
\emph{t-norm}) is an associative binary operation on $[0, 1]$
(henceforth, $I$) that is commutative, nondecreasing in each
place, such that $T(a, 1)=a$ for all $a\in I$.}
\end{definit}

\begin{definit}$\left.\right.$\vspace{.5pc} 

\noindent {\rm Let $T$ be a binary operation on $I$. Denote by
$T^*$ the function defined by $T^*(a,b):=1-T(1-a, 1-b)$ for all
$a, b\in I$. If  $T$ is a $t$-norm, then $T^*$ will be called the
\emph{t-conorm} of $T$. A~function $S$ is a $t$-conorm if there is
a $t$-norm $T$ such that $S=T^*$.}
\end{definit}
Clearly, $T^*$ is itself a binary operation on $I$, and
$T^{**}=T$. Instances of such $t$-norms and $t$-conorms are $M$
and $M^*$, respectively, defined by $M(x,y)=\min(x,y)$ and
$M^*(x,y)=\max(x,y)$.

\begin{definit}$\left.\right.$\vspace{.5pc} 

\noindent {\rm A \emph{triangle function} $\tau$ is an associative
binary operation on $\Delta^+$ that is commutative, nondecreasing
in each place, and has $\varepsilon_0$ as identity.}
\end{definit}
Also we let $\tau^1=\tau$ and
\begin{equation*}
\tau^n(F_1,\dots,F_{n+1})=\tau(\tau^{n-1}(F_1,\dots,F_n),F_{n+1})\;\;
\text{for}\; n\geq 2.
\end{equation*}

Let $T$ be a left-continuous $t$-norm and $T^*$ a right-continuous
$t$-conorm. Then instances of such triangle functions are $\tau_T$
and $\tau_{T^*}$ defined for all $F$, $ G\in\Delta^+ $ and every
$x\in \mathbf R^+$, respectively, by
\begin{equation*}
\tau_T(F,G)(x)=\sup\{T(F(u),G(v))\mid u+v=x\}
\end{equation*}
and
\begin{equation*}
\tau_{T^*}(F,G)(x)=\ell^-\,\inf \{T^*(F(u),G(v))\mid u+v=x\}.
\end{equation*}
The triangular function $\tau$ is said to be \emph{Archimedean} on
$\Delta^+$ if $\tau(F,G)<F$ for any $F$, $G$ in $\Delta^+$, such
that $F\neq \varepsilon_\infty$ and $G\neq \varepsilon_0$.

\begin{definit}$\left.\right.$\vspace{.5pc}  

\noindent {\rm Let $\tau_{1},\tau_{2}$ be two triangle functions.
Then $\tau_{1}$ dominates $\tau_{2}$, and we write $\tau_{1}\gg
\tau_{2}$, if for all $F_{1},F_{2},G_{1},G_{2}\in \Delta^{+}$,
\begin{equation*}
\tau_{1}(\tau_{2}(F_{1},G_{1}),\tau_{2}(F_{2},G_{2}))\geq
\tau_{2}(\tau_{1}(F_{1},F_{2}),\tau_{1}(G_{1},G_{2})).
\end{equation*}}
\end{definit}
In 1993, Alsina, Schweizer and Sklar \cite{alsina1} gave a new
definition of a probabilistic normed space as follows:

\begin{definit}$\left.\right.$\vspace{.5pc} 

\noindent {\rm A \emph{probabilistic normed space}, briefly a PN
space, is a quadruple $(V,\nu,\tau,\tau^*)$ in which $V$ is a
linear space, $\tau$ and $\tau^*$ are continuous triangle
functions with $\tau\leq\tau^*$ and $\nu$, the probabilistic norm,
is a map $\nu\hbox{\rm :}\ V\to\Delta^{+}$ such that
\begin{enumerate}
\renewcommand\labelenumi{(N\arabic{enumi})}
\leftskip 1.7pc
\item $\nu_p=\varepsilon_0$ if and only if $p=\theta$,
$\theta$ being the null vector in $V$;

\item $\nu_{-p}=\nu_p$\quad for every $p\in V$;

\item $\nu_{p+q}\geq\tau(\nu_p,\nu_q)$ for all $p,q \in V$;

\item $\nu_p\leq\tau^*(\nu_{\alpha p},\nu_{(1-\alpha)p})$ for every
$\alpha\in [0,1]$ and for every $p\in V$.\vspace{-.7pc}
\end{enumerate}}
\end{definit}

If, instead of (N1), we only have $\nu_{\theta}=\varepsilon_0$,
then we shall speak of a {\it probabilistic pseudo normed space},
briefly a PPN space. If the inequality (N4) is replaced by the
equality $\nu_p=\tau_M(\nu_{\alpha p},\nu_{(1-\alpha)p})$, then
the PN space is called a \emph{\v{S}erstnev space}; in this case,
a condition stronger than (N2)\,holds, namely
\begin{equation*}
\nu_{\lambda p}=\nu_p\left(\frac{j}{|\lambda |}\right), \qquad
\forall\lambda\neq 0, \:\forall p\in V;
\end{equation*}
here $j$ is the identity map on $\mathbf{R}$. A~\v{S}erstnev space
is denoted by $(V,\nu,\tau)$.

There is a natural topology in a PN space $(V,\nu,\tau,\tau^*)$,
called the {\it strong topology}; it is defined, for $t>0$, by the
neighbourhoods
\begin{equation*}
N_p(t):=\{q\in V\hbox{\rm :}\ d_S(\nu_{q-p},\epsilon_0)<t\}=\{q\in
V\hbox{\rm :}\ \nu_{q-p}(t)>1-t\}.
\end{equation*}

The strong neighbourhood system for $V$ is the union
$\bigcup_{p\in V}\mathcal{N }_p(\lambda)$ where $\mathcal{N}_p =
\{N_p(\lambda)\hbox{\rm :}\ \lambda>0\}$. The strong neighborhood
system for $V$ determines a Hausdorff topology for $V$.

A linear map $T\hbox{\rm :}\ (V,\nu,\tau,\tau^*)\rightarrow
(V',\nu',\sigma,\sigma^*)$, is said to be \emph{strongly bounded},
if there exists a constant $k>0$ such that, for all $p\in V $ and
$x>0$,
\begin{equation*}
\nu'_{Tp}(x)\geq \nu_p(x/k).
\end{equation*}

\begin{definit}$\left.\right.$\vspace{.5pc} 

\noindent {\rm A \emph{Menger PN space} is a PN space $(V,\nu,
\tau,\tau^*)$ in which $\tau=\tau_T$ and $\tau^*=\tau_{T^*}$ for
some $t$-norm $T$ and its $t$-conorm $T^*$. It will be denoted by
$(V,\nu,T)$.}
\end{definit}

\begin{definit}$\left.\right.$\vspace{.5pc} 

\noindent {\rm Let $(V,\nu,\tau,\tau^{*})$ be a PN space.
A~sequence $(p_{n})_{n}$ in $V$ is said to be \emph{strongly
convergent} to $p$ in $V$ if for each $\lambda>0$, there exists a
positive integer $N$ such that $p_{n}\in N_{p}(\lambda )$, for
$n\geq N$. Also the sequence $(p_{n})_{n}$ in $V$ is called a
\emph{strong Cauchy sequence} if, for every $\lambda >0$, there is
a positive integer $N$ such that $\nu_{p_{n}-p_{m}}(\lambda)>
1-\lambda$, whenever $m,n>N$. A~PN space $(V,\nu,\tau,\tau^{*})$
is said to be \emph{strongly complete} in the strong topology if
and only if every strong Cauchy sequence in $V$ is strongly
convergent  to a point in $V$.}
\end{definit}

\begin{lem}\hskip -.3pc \hbox{\rm \cite{alsina2}.}\ \ \label{alpha} 
If $|\alpha|\leq |\beta|,$ then $\nu_{\beta p}\leq \nu_{\alpha p}$
for every $p$ in $V$.
\end{lem}

\begin{definit}$\left.\right.$\vspace{.5pc}

\noindent {\rm Let $(V_1,\nu_1,\tau,\tau^{*})$ and $(V_2,\nu_2,
\tau,\tau^{*})$ be two PN spaces under the same triangle functions
$\tau$ and $\tau^*$. Let $\sigma$ be a triangle function. The
$\sigma$-product of the two PN spaces is the quadruple
\begin{equation*}
(V_1\times V_2, \nu_1 \sigma \nu_2,\tau,\tau^*),
\end{equation*}
where
\begin{equation*}
\nu_1 \sigma \nu_2\hbox{\rm :}\ V_1\times V_2\longrightarrow
\Delta^+
\end{equation*}
is a probabilistic semi-norm given by
\begin{equation*}
(\nu_1 \sigma \nu_2)(p,q):=\sigma(\nu_1(p),\nu_2(q))
\end{equation*}
for all $(p,q)\in V_1\times V_2$.}
\end{definit}

\section{Quotient PN space}\label{S:quot}

\setcounter{theore}{0}

According to \cite{schw-sk1} (see Definition~12.9.3 in p.~215),
one has the following:

\begin{definit}$\left.\right.$\vspace{.5pc} 

\noindent {\rm A triangle function $\tau$ is \emph{sup-continuous}
if, for every family $\{F_{\lambda}\hbox{\rm :}\ \lambda \in
\Lambda\}$ of \dfe 's in $\Delta^+$ and every $G\in \Delta^+$, }
\begin{equation*}
\sup_{\lambda \in \Lambda}\tau(F_\lambda, G)=\tau
\left(\sup_{\lambda \in \Lambda} F_\lambda,G\right).
\end{equation*}
\end{definit}
In view of Lemma~4.3.5 of \cite{schw-sk1}, this supremum is in
$\Delta^+$. An example of a \textit{sup-continuous} triangle
function is $\tau_T$, where $T$ is a left continuous $t$-norm.

\begin{definit}$\left.\right.$\vspace{.5pc} 

\noindent {\rm Let $W$ be a linear subspace of $V$ and denote by
$\sim_W$ a relation on the set $V$ defined via
\begin{equation*}
p_{\sim_W} q \Leftrightarrow p-q \in W,
\end{equation*}
for every $p,q \in V$.}
\end{definit}
Obviously this relationship is an equivalence relation and
therefore the set $V$ is partitioned into equivalence classes,
$V/\!\sim_{W}$.

\begin{propo}$\left.\right.$\vspace{.5pc} 

\noindent Let $(V,\nu,\tau,\tau^*)$ be a PN space. Suppose that
$\tau$ and $\tau^*$ are sup-continuous. Let $W$ be a subspace of
$V$ and \hbox{$V/\!\sim_W$} its quotient defined by means of the
equivalence relation $\sim_W$. Let $\nu'$ be the restriction of
$\nu$ to $W$ and define the mapping \hbox{$\bar{\nu}\hbox{\rm :}\
V/\!\sim_W \,\, \rightarrow \Delta^+,$} for all $p\in V,$ by
\begin{equation*}
\bar{\nu}_{p+W}(x):= \sup_{q \in W}\{\nu_{p+q}(x)\}.
\end{equation*}
Then{\rm ,} $(W,\nu',\tau,\tau^*)$ is a PN space and
$(V/\!\sim_W,\bar{\nu},\tau,\tau^*)$ is a PPN space.
\end{propo}

\begin{proof}
The first statement is immediate. The remainder of the theorem is
guaranteed by the fact that $W$ is not necessarily closed in the
strong topology.\hfill $\Box$
\end{proof}
Notice that by Lemma~4.3.5 of \cite{schw-sk1}, $\bar{\nu}_{p+W}$
is in $\Delta^+$.

Hereafter we denote by $p_W$ the subset $p+W$ of $V$, \ie\ an
element of quotient, and the strong neighbourhood of $ p_W$ by
$N'_{p_W}(t)$.

\begin{theor}[\!] 
Let $W$ be a linear subspace of $V$. Then the following statements
are equivalent{\rm :}\vspace{-.3pc}
\begin{enumerate}
\renewcommand\labelenumi{\rm (\alph{enumi})}
\leftskip .1pc
\item  $(V/\!\sim_W, \bar{ \nu },\tau,\tau^*)$ is a PN
space{\rm ;}

\item $W$ is closed in the strong topology of
$(V,\nu,\tau,\tau^*)$.\vspace{-.5pc}
\end{enumerate}
\end{theor}

\begin{proof}
Let $(V,\nu,\tau,\tau^*)$ be a PN space. For every $p$ in the
closure of $W$ and for each $n\in \Na$ choose $q_n\in N_p(1/n)\cap
W$. Then
\begin{equation*}
\bar{ \nu }_{p_W}(1/n)= \sup_{q\in W}\nu_{p+q}(1/n)\geq
\nu_{p-q_n}(1/n)>1-1/n,
\end{equation*}
and hence, $d_S(\bar{ \nu }_{p_W},\varepsilon_0)<1/n$. Thus
$p_W=W$ and hence, $p\in W$ and $W$ is closed.

Conversely, if $W$ is closed, let $p\in V$ be such that
$\bar{\nu}_{p_W}=\varepsilon_0.$ If $p \not \in W$, then
$N_p(t)\cap W=\varnothing$, for some $t>0$. That is to say, for
every $q\in W$, $\nu_{p-q}(t)\leq 1-t$. Therefore
$\bar{\nu}_{p_W}(t)=\sup_{q\in W}\nu_{p+q}(t)\leq 1-t$, which is a
contradiction.\hfill $\Box$
\end{proof}

It is of interest to know  whether a PN space can be obtained from
a PPN space. An affirmative answer is provided by the following
proposition.

\begin{propo}$\left.\right.$\vspace{.5pc} 

\noindent Let $(V,\nu,\tau,\tau^*)$ be a PPN space and define
\begin{equation*}
C=\{p\in V\hbox{\rm :}\ \nu_p= \varepsilon_0\}.
\end{equation*}
Then $C$ is the smallest closed subspace of $(V,\nu,\tau,\tau^*)$.
\end{propo}

\begin{proof}
If $p,q\in C$, then $p+q \in C$ because $\nu_{p+q}\geq
\tau(\nu_p,\nu_q)=\varepsilon_0$. Now suppose $p \in C$. For
$\alpha \in [0,1]$ one has $\nu_{\alpha p}\geq \nu_p$ by
Lemma~\ref{alpha}. For $\alpha >1$, let $k=[\alpha] +1$. Then,
using the iterates of (N3) one has, $\nu_{kp}\geq
\tau^{k-1}(\nu_p,\dots,\nu_p)=\varepsilon_0$. By the
above-mentioned lemma one has $\nu_{\alpha p}\geq \nu_{kp}$. As a
consequence, $\alpha p$ belongs to $C$ for all $\alpha \in
\mathbf{R}$.

Furthermore it is easy to check that the set $C$ is closed because
of the continuity of the probabilistic norm, $\nu$ (see Theorem~1
in \cite{alsina2}).

Now, let $W$ be a closed linear subspace of $V$ and $p\in C$.
Suppose that for some $t>0$, $N_p(t)\cap W=\varnothing$, then
$\nu_{p}(t)\leq 1-t$, which is a contradiction; hence $C\subseteq
W$.\hfill $\Box$
\end{proof}

\begin{remar} 
{\rm Moreover, with $V$ and $C$ as in Proposition~3.5, for all $p
\in V$ and $r \in C$, one has
\begin{equation*}
\bar{ \nu }_{p_W}\geq
\nu_p=\nu_{p+r-r}\geq\tau(\nu_{p+r},\nu_{-r})=\nu_{p+r}.
\end{equation*}
Thus the probabilistic norm $\bar{ \nu }$ in $(V/\!\sim_{C},
\bar{\nu},\tau, \tau^*)$ coincides with that of $(V, \nu, \tau,
\tau^*)$.}
\end{remar}

\begin{exampl} 
{\rm Let $(V,\nu,T)$ be a Menger PN space. Suppose that $ W$ is a
closed subspace of $V$, and \hbox{$V/\!\sim_W$} its quotient. Then
$(W,\nu',T)$ and \hbox{$(V/\!\sim_W, \bar{\nu}, T)$} are Menger
PN\break spaces.}
\end{exampl}

\begin{coro}$\left.\right.$\vspace{.5pc} 

\noindent Let $(V,\nu,\tau,\tau^*)$ be a \v{S}erstnev PN space.
Suppose that $\tau$ is sup-continuous. Let $W$ be a closed
subspace of $V$ and \hbox{$V/\!\sim_W$} its quotient. Then{\rm ,}
$(W,\nu',\tau,\tau^*)$ and
\hbox{$(V/\!\sim_W,\bar{\nu},\tau,\tau^*)$} are \v{S}erstnev PN
spaces.
\end{coro}

\begin{theor}[\!] 
Let $(V,\nu,\tau,\tau^*)$ be a PN space. Suppose that $\tau$ and
$\tau^*$ are sup-continuous. Let $W$ be a closed subspace of $V$
with respect to the strong topology of $(V,\nu,\tau,\tau^*)$. Let
\begin{equation*}
\pi\hbox{\rm :}\ V\rightarrow V/\!\sim_W
\end{equation*}
be the canonical projection. Then $\pi$ is strongly bounded{\rm ,}
open{\rm ,} and continuous with respect to the strong topologies
of $(V,\nu,\tau,\tau^*)$ and
\hbox{$(V/\!\sim_W,\bar{\nu},\tau,\tau^*)$}. In addition{\rm ,}
the strong topology and the quotient topology on
\hbox{$V/\!\sim_W,$} induced by $\pi,$ coincide.
\end{theor}

\begin{proof}
One has that $\bar{\nu}_{p_W}\geq \nu_p$ which implies $\pi$ is
strongly bounded, and hence continuous (see Theorem~3.3 in
\cite{LRS3}).

The map $\pi$ is open because of the equality
$\pi(N_p(t))=N'_{p_W}(t)$.\hfill $\Box$
\end{proof}

\begin{exampl} 
{\rm Let $(V,\norm{\cdot})$ be a normed space and define
$\nu\hbox{:}\ V\rightarrow \Delta^+$ via
$\nu_p:=\varepsilon_{\parallel p\parallel}$ for every $p \in V$.
Let $\tau,\tau^*$ be continuous triangle functions such that
$\tau\leq\tau^*$ and $\tau(\varepsilon_a, \varepsilon_b)=
\varepsilon_{a+b}$, for all $a,b>0$. For instance, it suffices to
take $\tau=\tau_T$ and $\tau^* =\tau_{T^*}$, where $T$ is a
continuous $t$-norm and $T^*$ is its $t$-conorm. Then $(V, \nu,
\tau,\tau^*)$ is a PN space (see Example~1.1 of \cite{LRS3}).

Assume that $\tau$ is sup-continuous. Let W be a closed linear
subspace of $V$ with respect to the strong topology of $(V, \nu,
\tau, \tau^*)$. By Theorem~3.4, \hbox{$(V/\!\sim_W, \mu, \tau,
\tau^*)$} is a PN space in which, $\mu_{p_W}=\sup_{w\in W}
\varepsilon_{\parallel p+w\parallel}$. On the other hand, if one
considers the normed space \hbox{$(V/\!\sim_W,\norm{\cdot}')$},
where $\norm{p_W}'=\inf_{w\in W}\norm{p+w}$, then one can easily
prove that the PN structure given to the normed space
\hbox{$(V/\!\sim_W, \norm{\cdot}')$} by means of
$\eta_{p_W}:=\varepsilon_{\parallel p_W \parallel'}$ coincides
with $(V/\!\sim_W, \eta, \tau,\tau^*)$.}
\end{exampl}

\section{Completeness results}\label{S:compl}

\setcounter{theore}{0}

Here we study the completeness of a quotient PN space. When a PN
space $(V, \nu, \tau, \tau^*)$ is strongly complete, then we say
that it is a probabilistic normed Banach (henceforth PNB) space.

\begin{lem} 
Given the PN space \hbox{$(V/\!\sim_W,\bar{\nu},\tau,\tau^*)$} in
which $\tau$ and $\tau^*$ are sup-continuous{\rm ,} let $W$ be a
closed subspace of $V$.
\begin{enumerate}
\renewcommand\labelenumi{\rm (\roman{enumi})}
\leftskip .15pc
\item If $p\in V,$ then for every $\epsilon >0$ there is a $p'$ in
$V$ such that $p'+W=p+W$ and
\begin{equation*}
\hskip -1.25pc d_S(\nu_{p'},\varepsilon_0) <
d_S(\bar{\nu}_{p+W},\varepsilon_0) +\epsilon.
\end{equation*}

\item If $p$ is in $V$ and $ \bar{\nu}_{p+W}\geq G $ for some d.f.
$G \neq \varepsilon_0,$ then there exists $p'\in V$ such that
$p+W=p'+W$ and $\nu_{p'}\geq\tau(\bar{\nu}_{p+W},G)$.
\end{enumerate}
\end{lem}

\begin{proof}$\left.\right.$
\begin{enumerate}
\renewcommand\labelenumi{\rm (\roman{enumi})}
\leftskip .15pc
\item We know
\begin{equation*}
\hskip -1.25pc \bar{\nu}_{p+W}=\sup\{\nu_{p-q}\hbox{\rm :}\ q\in
W\}.
\end{equation*}
Now, let $q$ be an element of $W$ such that
\begin{equation*}
\hskip -1.25pc \bar{\nu}_{p+W}< \nu_{p-q} + \frac{\epsilon}{2}.
\end{equation*}
We put $p-q=p'$. Now,
\begin{align*}
\hskip -1.25pc d_S(\bar{\nu}_{p+W},\varepsilon_0) &=
\inf\{h>0\hbox{:}\
\bar{\nu}_{p+W}(h^+)> 1-h\}\\[.3pc]
&\geq \inf \left\{h>0\hbox{:}\ \nu_{p'}(h^+)+\frac{\epsilon}{2} >
1-h\right\}\\[.3pc]
&= \inf \left\{h>0\hbox{:}\ \nu_{p'}(h^+) > 1-\left( h+\frac{\epsilon}{2}\right)\right\}\\[.3pc]
&\geq \inf \left\{h>0\hbox{:}\ \nu_{p'}\left( \left( h +
\frac{\epsilon}{2} \right)^+\right) > 1- \left(h+\frac{\epsilon}{2}\right)\right\}\\[.3pc]
&> d_S(\nu_{p'},\varepsilon_0) - \epsilon.
\end{align*}

\item Because of the definition of supremum and sup-continuity of
$\tau $, there exists a $q_n\in W$ such that $q_n\rightarrow q$ if
$n\rightarrow +\infty$ and
\begin{equation*}
\nu_{p+q_n}>\tau(\bar{\nu}_{p_W},\varepsilon_0)-\frac{1}{n}\geq
\tau(\bar{\nu}_{p_W},G)-\frac{1}{n}.
\end{equation*}

Now it is enough to put $p'=p+q$ and see that, when $n\rightarrow
+\infty$, one has $\nu_{p+q} \geq \tau(\bar{\nu}_{p_W},G)$.\hfill
$\Box$\vspace{-2.3pc}
\end{enumerate}
\end{proof}
\pagebreak

Let $p,q$ be elements of $V$ such that $d_S(\nu_{(p-q)+W},
\varepsilon_0) < \delta$ for some positive $\delta$. By Lemma~4.1,
there is a $q'\in V$ such that $(p-q')+W=(p-q)+W$ and
\begin{align*}
d_S(\nu_{p-q'},\varepsilon_0) < \delta.
\end{align*}

\begin{theor}[\!] 
Let $W$ be a closed subspace of $V$ and suppose that $(V, \nu,
\tau, \tau^*)$ is a PNB space with $\tau$ and $\tau^*$
sup-continuous. Then{\rm ,} $(V/\!\sim_W, \bar{\nu}, \tau,
\tau^*)$ is also a PNB space.
\end{theor}

\begin{proof}
Let $(a_n)$ be a strong Cauchy sequence in \hbox{$(V/\!\sim_W,
\bar{\nu}, \tau, \tau^*)$}, \ie\ for every $\delta >0$, there
exists $n_0=n_0(\delta) \in \mathbf N$ such that, for all $m,n>
n_0$,
\begin{equation*}
d_S(\bar{\nu}_{a_n-a_m},\varepsilon_0)<\delta.
\end{equation*}

Now, define a strictly  decreasing sequence $(\delta_n)$ with
$\delta_n > 0$ in the following way: let $\delta _1>0$ be such
that $\tau(B_{d_S}(\varepsilon_0;\delta_1)\times B_{d_S}
(\varepsilon_0;\delta_1))\subseteq B_{d_S}(\varepsilon_0;1 )$
where $B_{d_S}(\varepsilon_0;\lambda )=\{F\in
\Delta^+;d_S(F,\varepsilon_0)<\lambda\}$. For $n>1$, define
$\delta_n$ by induction in such a manner that
\begin{equation}
\tau(B_{d_S}(\varepsilon_0;\delta_n)\times B_{d_S} (\varepsilon_0;
\delta_n))\subseteq B_{d_S}\left(\varepsilon_0; \min\left(
\frac{1}{n}, \delta_{n-1}\right)\right).
\end{equation}
There is a subsequence $(a_{n_i})$ of $(a_n)$ with
\begin{equation}
d_S(\bar{\nu}_{a_{n_i}-a_{n_{i+1}}},\varepsilon_0)<\delta_{i+1}.
\end{equation}
Because of the definition of the canonical projection $\pi$ one
can say that $\pi^{-1}(N'_{p_W}(t))=N_p(t)$ and consequently
$\pi^{-1}(a_{n_i})=x_i$ exists. Inductively, from Lemma~4.1 we can
find $x_i\in V$ such that $\pi(x_i)=a_{n_i}$ and then
\begin{equation}
d_S(\nu_{x_i-x_{i+1}},\varepsilon_0)<\delta_{i+1}
\end{equation}
holds. We claim that $(x_i)$ is a strong Cauchy sequence in
$(V,\nu, \tau, \tau^*)$. By applying the relations (1), (2) and
(3) to $i = m-1$ and $i=m-2$, and using Lemma~4.3.4 of
\cite{schw-sk1}, one obtains the inequalities
\begin{align*}
d_S(\nu_{x_m-x_{m-2}},\varepsilon_0) &\leq d_S (\tau(
\nu_{x_{m-1}-x_m},\nu_{x_{m-2}-x_{m-1}}),\varepsilon_0)\\[.3pc]
&< \min\left(\frac{1}{m-1},\delta_{m-2}\right).
\end{align*}
Following this reasoning, we obtain that $d_S(\nu_{x_m-x_n},
\varepsilon_0)<1/n$ and therefore, $(x_i)$ is a strong Cauchy
sequence. Since it was assumed that $(V, \nu, \tau, \tau^*)$ is
strongly complete, $(x_i)$ is strongly convergent and hence, by
the continuity of $\pi$, $(a_{n_i})$ is also strongly convergent.
From this and taking into account the continuity of $\tau$ and
Lemma~4.3.4 of \cite{schw-sk1}, one sees that the whole sequence
$(a_n)$ strongly converges.\hfill $\Box$
\end{proof}

The converse of the above theorem also holds.

\begin{theor}[\!] 
Let $(V,\nu,\tau,\tau^{*})$ be a PN space in which $\tau$ and
$\tau^*$ sup-continuous{\rm ,} and let \hbox{$(V/\!\sim_W,
\bar{\nu}, \tau,\tau^{*})$} be its quotient space with respect to
the closed subspace $W$. If any two of the three spaces $V,$ $W$
and $V/\!\sim_W$ are strongly complete{\rm ,} so is the third.
\end{theor}

\begin{proof}
If $V$ is a strongly complete PN space, so are \hbox{$V/\!\sim_W$}
and $W$. Therefore all one needs to check is that $V$ is strongly
complete whenever both $W$ and \hbox{$V/\!\sim_W$} are strongly
complete. Suppose $W$ and \hbox{$V/\!\sim_W$} are strongly
complete PN spaces and $(p_{n})$ be a strong Cauchy sequence in
$V$. Since
\begin{equation*}
\bar{\nu}_{(p_{m}-p_{n})+W}\geq \nu_{p_{m}-p_{n}}
\end{equation*}
whenever $m,n\in {\bf N}$, the sequence $(p_{n}+W)$ is strong
Cauchy in \hbox{$V/\!\sim_W$} and, therefore, it strongly
converges to $q+W$ for some $q\in V$. Thus there exists a sequence
of d.f's $(H_{n})$ such that $H_{n}\longrightarrow
\varepsilon_{0}$ and $\bar{\nu}_{(p_{n}-q)+W}>H_{n}$. Now by
Lemma~4.1 there exists $(q_{n})$ in $V$ such that
$q_{n}+W=(p_{n}-q)+W$ and
\begin{equation*}
\nu_{q_{n}}>\tau(\bar{\nu}_{(p_{n}-q)+W},H_{n}).
\end{equation*}
Thus $\nu_{q_{n}}\longrightarrow \varepsilon_{0}$ and consequently
$q_{n}\longrightarrow\theta$. Therefore $(p_{n}-q_{n}-q)$ is a
strong Cauchy sequence in $W$ and is strongly convergent to a
point $r\in W$ and implies that $(p_{n})$ strongly converges to
$r+q$ in $V$. Hence $V$ is strongly complete.\hfill $\Box$
\end{proof}

\begin{theor}[\!] 
Let $(V_1,\nu^1,\tau,\tau^*), \dots, (V_n,\nu^n,\tau,\tau^*)$ be
PNB spaces in which $\tau$ and $\tau^*$ are sup-continuous.
Suppose that there is a triangle function $\sigma$ such that
$\tau^*\gg\sigma$ and $\sigma\gg\tau$. Then their $\sigma$-product
is a PNB space.
\end{theor}

\begin{proof} One proves for $n=2$ (see Theorem~2 in \cite{L}), and
then we apply induction for an arbitrary $n$. Since the quotient
norm of
\begin{equation*}
\frac{V_1\times V_2}{V_1\times{\theta_2}}\ (\simeq V_2)
\end{equation*}
is the same as $\nu^2$ and the restriction of the product norm of
$V_1 \times V_2$ to $V_1\times {\theta_2}(\simeq\!V_1)$ is the
same as $\nu^1$ (see \cite{L}), and in view of Theorem~4.3, the
proof is complete.\hfill $\Box$
\end{proof}

By Theorem~3.9 the following corollaries can be proved easily.

\begin{coro}$\left.\right.$\vspace{.5pc} 

\noindent Under the assumptions of Proposition~{\rm 3.3} and if
$W$ is a closed subset of $V,$ the probabilistic norm
$\bar{\nu}\hbox{\rm :}\ V/\!\sim_W \, \rightarrow \Delta^+$ in
$(V/\!\sim_W, \bar{\nu}, \tau, \tau^*)$ is uniformly continuous.
\end{coro}

\begin{proof}
Let $\eta$ be a positive real number, $\eta>0$. By Theorem~3.9
there exists a pair $(p', q')$ in $(V \times V)$ such that
$d_S(\bar{\nu}_{\pi(p-p')},\varepsilon_{0})<\eta$ and
$d_S(\bar{\nu}_{\pi(q-q')},\varepsilon_{0})<\eta$, whenever
$d_S(\nu_{p-p'},\varepsilon_{0})<\eta$ and $d_S(\nu_{q-q'},
\varepsilon_{0})<\eta$.

On the other hand, we have
\begin{equation*}
\bar{\nu}_{\pi(p'-q')}\geq \tau(\tau(\bar{\nu}_{\pi(p-p')},
\bar{\nu}_{\pi(q-q')}),\bar{\nu}_{\pi(p-q)})
\end{equation*}
and
\begin{equation*}
\bar{\nu}_{\pi(p-q)}\geq \tau(\tau(\bar{\nu}_{\pi(p-p')},
\bar{\nu}_{\pi(q-q')}),\bar{\nu}_{\pi(p'-q')}).
\end{equation*}
Thus, from the relationship (12.1.5) and Lemma~12.2.1 in
\cite{schw-sk1} it follows that for any $h>0$ there is an
appropriate $t > 0$ such that
\begin{equation*}
d_S(\bar{\nu}_{\pi(p-q)},\bar{\nu}_{\pi(p'-q')})<h,
\end{equation*}

$\left.\right.$\vspace{-1.3pc}

\pagebreak

\noindent whenever $p'\in N_p(\eta)$ and $q'\in N_q(\eta)$. This
implies that $\bar{\nu}$ is a uniformly continuous mapping from
$V/\!\sim_W$ into $\Delta^{+}$.\hfill $\Box$
\end{proof}

Also the inequality $d_S(\bar{\nu}_{\pi((p+q)- (p'+q'))},
\varepsilon_{0}) \leq d_S(\nu_{(p+q)-(p'+q')},\varepsilon_{0})$
implies that $(V/\!\sim_W,+)$ is a topological group.

\begin{coro}$\left.\right.$\vspace{.5pc} 

\noindent Let $(V,\nu,\tau,\tau^{*})$ be a PN space such that
$\tau^{*}$ is Archimedean{\rm ,} $\tau $ and $\tau^*$ are
sup-continuous{\rm ,} and $\nu_{p}\neq \varepsilon_{\infty}$ for
all $p\in V$. If we define quotient probabilistic norm via
Proposition~{\rm 3.3,} then
\hbox{$(V/\!\sim_W,\bar{\nu},\tau,\tau^{*})$} is a PPN space where
the scalar multiplication is a continuous mapping from
\hbox{$R\times V/\!\sim_W$} into $V/\!\sim_W$.
\end{coro}

\begin{proof}
For any $p\in V$ and $\alpha ,\beta\in R$ we know
$d_S(\bar{\nu}_{\pi(\alpha p)},\nu_{\pi(\beta p)})$ is small
whenever $d_S(\bar{\nu}_{\pi((\alpha -\beta) p)},\varepsilon_{0})$
is small. But
\begin{equation*}
d_S(\bar{\nu}_{\pi((\alpha -\beta) p)},\varepsilon_{0})\leq
d_S(\nu_{(\alpha -\beta) p},\varepsilon_{0})
\end{equation*}
and by Lemma~3 of \cite{alsina2}, $d_S(\nu_{(\alpha -\beta)
p},\varepsilon_{0})$ is small whenever $|\alpha -\beta|$ is
small.\hfill $\Box$
\end{proof}

\section*{Acknowledgments}

The authors would like to thank the referee for giving useful
comments and suggestions for the improvement of this paper. The
authors wish to thank C~Sempi for his helpful suggestions. The
first author was supported by grants from the Ministerio de
Ciencia y Tecnolog\'{\i}a (BFM2005-06522) and from the Junta de
Andaluc\'{\i}a (CEC-JA FQM-197), Spain.

\end{document}